\documentclass[12pt,reqno]{amsart}

\usepackage[margin=1in]{geometry}  
\usepackage{graphicx}              
\usepackage{amsmath}               
\usepackage{amsfonts}              
\usepackage{amsthm}                
\usepackage{amssymb}
\usepackage{appendix}
\usepackage[utf8]{inputenc}
\usepackage{mathabx}
\usepackage{floatrow}

\newtheorem{thm}{Theorem}[section]

\newtheorem{conj}[thm]{Conjecture}
\newtheorem{que}[thm]{Question}

\newcommand{\bd}[1]{\mathbf{#1}}  
\newcommand{\CC}{\mathbb{C}}

\newcommand{\al}[1]{\begin{align}#1\end{align}}

\begin{document}

\title{Spatial Statistics of Apollonian Gaskets}

\author{Weiru Chen}
\
\author{Mo Jiao}
\
\author{Calvin Kessler}
\
\author{Amita Malik}
\
\author{Xin Zhang}
\address{University of Illinois at Urbana-Champaign, Department of Mathematics\\
1409 W Green Street, Urbana, IL 61801} \email{wchen108@illinois.edu\\  \newline
  \indent \indent \indent \indent \indent \indent \indent qramir2@illinois.edu\\
\newline   \indent \indent \indent \indent \indent \indent \indent
mojiao2@illinois.edu\\ \newline
  \indent \indent \indent \indent \indent \indent \indent
 amalik10@illinois.edu \\ \newline
  \indent \indent \indent \indent \indent \indent \indent
 xz87@illinois.edu
}

\begin{abstract} 
Apollonian gaskets are formed by repeatedly filling the interstices between four mutually tangent circles with further tangent circles. We experimentally study the pair correlation, electrostatic energy, and nearest neighbor spacing of centers of circles from Apollonian gaskets. Even though the centers of these circles are not uniformly distributed in any `ambient' space, after proper normalization, all these statistics seem to exhibit some interesting limiting behaviors. 
\end{abstract}

\maketitle

\section{introduction}
Apollonian gaskets, named after the ancient Greek mathematician, Apollonius of Perga (200 BC), are fractal sets obtained by starting from three mutually tangent circles and iteratively inscribing new circles in the curvilinear triangular gaps. Over the last decade, there has been a resurgent interest in the study of Apollonian gaskets. Due to its rich mathematical structure, this topic has attracted attention of experts from various fields including number theory, homogeneous dynamics, group theory, and significant results have been obtained.  \par

\begin{figure}[h]
\begin{center}
\includegraphics[scale=0.3]{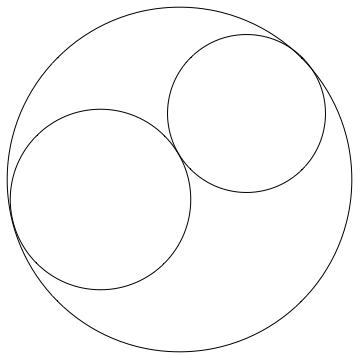}
\includegraphics[scale=0.3]{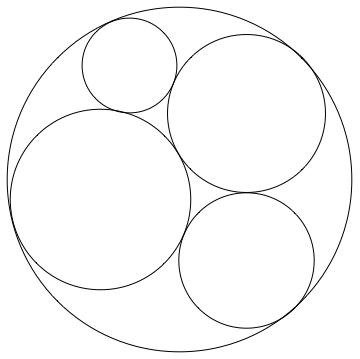}
\includegraphics[scale=0.3]{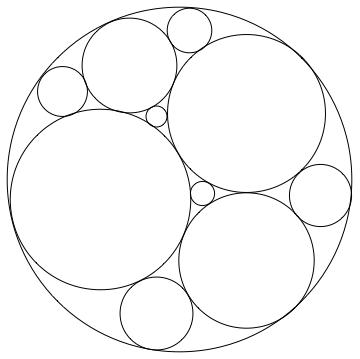}
\end{center}
\caption{Construction of an Apollonian gasket}
\label{fig:taxicab}
\end{figure}

\begin{figure}[!h]
           \begin{floatrow}
             \ffigbox{\includegraphics[scale = 0.22]{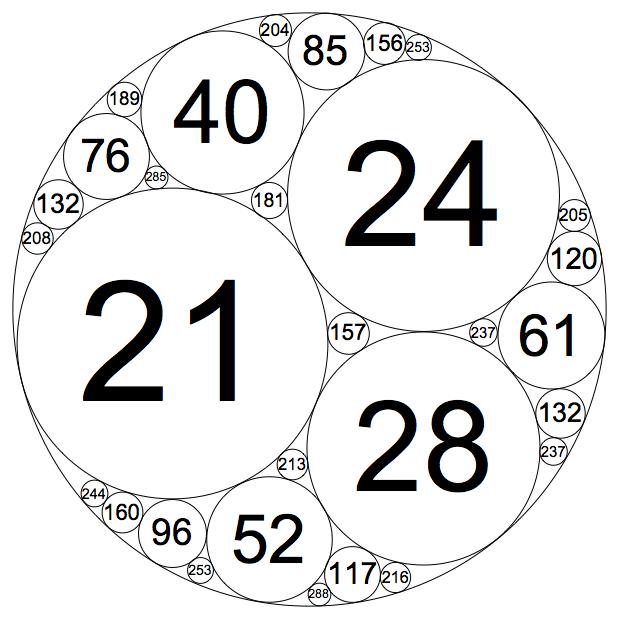}}{\caption{An integer Apollonian gasket}\label{apon}}
             \ffigbox{\includegraphics[scale = 0.18]{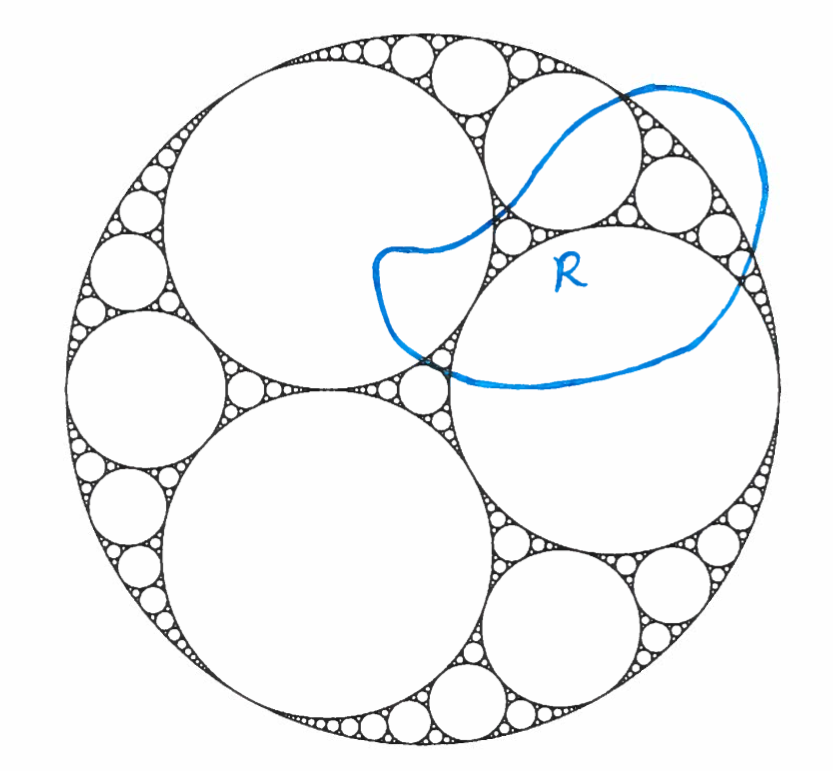}}{\caption{A region $\mathcal{R}$ with piecewise analytic boundary\label{aponr}}\label{center}}
           \end{floatrow}
        \end{figure}

For example, it has been known since Soddy \cite{So36} that there exist Apollonian gaskets with all circles having integer curvatures (reciprocal of radii). This is due to the fact that the curvatures from any four mutually tangent circles satisfy a quadratic equation (see Figure \ref{apon}). Inspired by \cite{GLMWY03}, \cite{Fu10}, and \cite{BF12}, Bourgain and Kontorovich used the circle method to prove a fascinating result that for any \emph{primitive} integral (integer curvatures with {\it{gcd} 1)} Apollonian gasket, almost every integer in certain congruence classes modulo 24 is a curvature of some circle in the gasket.  \par

In another direction, Kontorovich and Oh \cite{KO11} obtained an asymptotic result for counting circles from an Apollonian gasket $\mathcal{P}$ using the spectral theory of infinite volume hyperbolic spaces, which was originally developed by Lax-Phillips \cite{LP82}. Their result is stated below.
\begin{thm}[Kontorovich-Oh]
Fix an Apollonian gasket $\mathcal{P}$, and let $\mathcal P_T$ be the collection of circles with curvatures $<T$. Then as $T$ approaches infinity, 
$$\lim_{T\rightarrow\infty}\frac{\#\mathcal P_T}{T^{\delta}}= c_{\mathcal{P}},$$
where $c_{\mathcal P}$ is a positive constant depending on $\mathcal P$, and $\delta\approx 1.305688$ is the Hausdorff dimension of any Apollonian gasket.
\end{thm}
The reason for all Apollonian gaskets to have the same Hausdorff dimension is that they belong to the same conformally equivalent class. In other words, for any two fixed gaskets, one can always find a M\"obius transformation which takes one gasket to the other. The estimate $\delta\approx1.305688$ was obtained in \cite{Mc98}.\par
Kontorovich and Oh's result was refined by Oh and Shah \cite{OS12} using homogeneous dynamics. 
\begin{thm}[Oh-Shah]\label{0827}For any Apollonian gasket $\mathcal{P}$ placed in the complex plane $\mathbb C$, there exists a finite Borel measure $\mu$, such that for any region $\mathcal{R}\subset \mathbb C$ with piecewise analytic boundary (see Figure \ref{aponr}), the cardinality of $\mathcal P_T(\mathcal R)$, the set of circles from $\mathcal P_T$ lying in $\mathcal R$, satisfies 
$$\lim_{T\rightarrow \infty}\frac{\mathcal P_T(\mathcal R)}{T^{\delta}}=\mu(\mathcal R).$$
\end{thm}

Theorem \ref{0827} gives a satisfactory explanation on how circles are distributed in an Apollonian gasket in large scale. However, it yields little information on questions concerning the fine scale distribution of circles. For example, one such question is the following.
 
\begin{que}\label{que1}If one sits at the center of a random circle from $\mathcal{P}_T$, how many circles can one see within a distance of $10/T$?
\end{que}
 
The fine structures of Apollonian gaskets are encoded by local spatial statistics. In this article, we report our empirical results on some of such statistics, namely, pair correlation, nearest neighbor spacing and electrostatic energy. We find numerically that after proper normalization all these statistics exhibit some interesting limiting behavior when the growing parameter $T$ approaches infinity. Our conjectures in this direction are based on these numerical results. In particular, a proof of Conjecture \ref{con1} will provide an asymptotic answer to Question \ref{que1} when $T$ is large. \par
These spatial statistics have been widely used in disciplines such as physics and biology. For instance, in microscopic physics, the Kirkwood-Buff Solution Theory \cite{BK51} links the pair correlation function of gas molecules, which encodes the microscopic details of the distribution of these molecules, to macroscopic thermodynamical properties of the gas such as pressure and potential energy. On a macroscopic level, cosmologists use pair correlations to study the distribution of stars and galaxies.  \par
Within mathematics, there is a rich literature on the spatial statistics of point processes arising from various settings. A stunning application of pair correlation is a discovery made by Montgomery and Dyson. Montgomery computed that the pair correlation function of the zeros of the Riemann zeta function, which agrees with the pair correlation function of random Hermitian matrices computed by Dyson. This conjectural relation (still unproven) is known as the Montgomery's pair correlation conjecture \cite{Mo73}. It appears to give evidences to Hilbert and P\'olya's speculation that zeros of the Riemann zeta function correspond to eigenvalues of a self-adjoint operator on a Hilbert space.  \par
Spatial statistics from some other point processes have otherwise been rigorously established: gap distribution of the fractional parts of $(\sqrt{n})$ by Elkies and McMullen \cite{EM04}, distribution of directions of Euclidean or hyperbolic lattices \cite{BZ06}, \cite{BPZ14}, \cite{KK15}, \cite{RS14}, \cite{MV14},  distribution of Farey sequences \cite{Ha70}, \cite{BCZ01}, \cite{ABC01}, and gap distribution of saddle connection directions in translation surfaces \cite{AC10}, \cite{ACL15}. Our list of interesting works here is far from inclusive.  \par
There is a fundamental difference between all mentioned works above and our investigation of circles. In their work, the underlying point sequences become uniformly distributed in their `ambient' spaces.  In our case, we parametrize each circle by its center and we define the distance between two circles by the Euclidean distance of their centers, then our study of circles becomes the study of their centers. However, the set of centers is clearly not even dense in any reasonable ambient space such as $D$, the disk bounded by the largest circle of the gasket. In fact, we notice that centers tend to cluster over some tiny regions and meanwhile we can find in $D$ plenty of holes in which no center can be found. Consequently, we need different normalizations of parameters, as hinted in the last author's work \cite{Zhang15} on the gap distribution of a point orbit of an infinite-covolume Schottky group. \par

\section{Experimental results and conjectures}
The point process under our investigation is $\mathcal C_T$, the set of centers of circles with curvatures $<T$. In this section, we provide data for the (normalized) spatial statistics such as electrostatic energy, nearest neighbor spacing, pair correlations and state our conjectures. All packings under consideration here come from four mutually tangent circles $C_0,C_1,C_2,C_3$ with $C_0$ bounding the other three and of radius 1. We use $\mathbb C$-coordinates for these circles so that $C_0=\left\{z\in\mathbb C: \vert z\vert =1  \right\}$, and $C_1$ tangent to $C_0$ at $e^{0\bd{i}}=1$. Suppose $C_2$ and $C_3$ are tangent to $C_0$ at $e^{\theta_1\bd{i}}$ and $ e^{\theta_2\bd{i}}$, respectively. The pair $(\theta_1,\theta_2)$ then uniquely determines an Apollonian gasket which we denote by $\mathcal P(\theta_1,\theta_2)$.\par

\subsection{Pair correlation}
The pair correlation function $F_T(s)$ for the set $\mathcal C_T$ is defined to be 
$$F_T(s):=\frac{1}{2\#\mathcal C_T}\sum_{\substack{p,q\in\mathcal C_T\\p\neq q}}\bd{1}\left\{d(p,q)T< s\right\},$$
where $d(\cdot,\cdot)$ is the Euclidean distance function, and $\bd{1}\{\rm{A}\}$ is the indicator function which has the value 1 if the condition $A$ is true, and 0 otherwise. \par 
From Theorem \ref{0827}, we have $\#\mathcal{C}_T\sim c_{\mathcal P}T^{\delta}$. If these centers were randomly distributed in $D$, then a typical point has distance $\asymp 1/T^{{\delta}/{2}}$ to its nearest neighbor. But here we need to normalize the distance by multiplying $T$ instead of $T^{\delta/2}$. The reason is that in the family $\mathcal P_T$, a typical circle has radius $\asymp 1/T$, and the neighbor circles in $\mathcal{C}_T$ also typically have distance $\asymp 1/T$ (recall that the distance between two circles is the Euclidean distance between their centers), so that $1/T$ is the right scale to measure the local spacing of circles. We also use the same normalization for the nearest neighbor spacing statistics.\par

\begin{figure}[!h]
           \begin{floatrow}
             \ffigbox{\includegraphics[scale = 0.27]{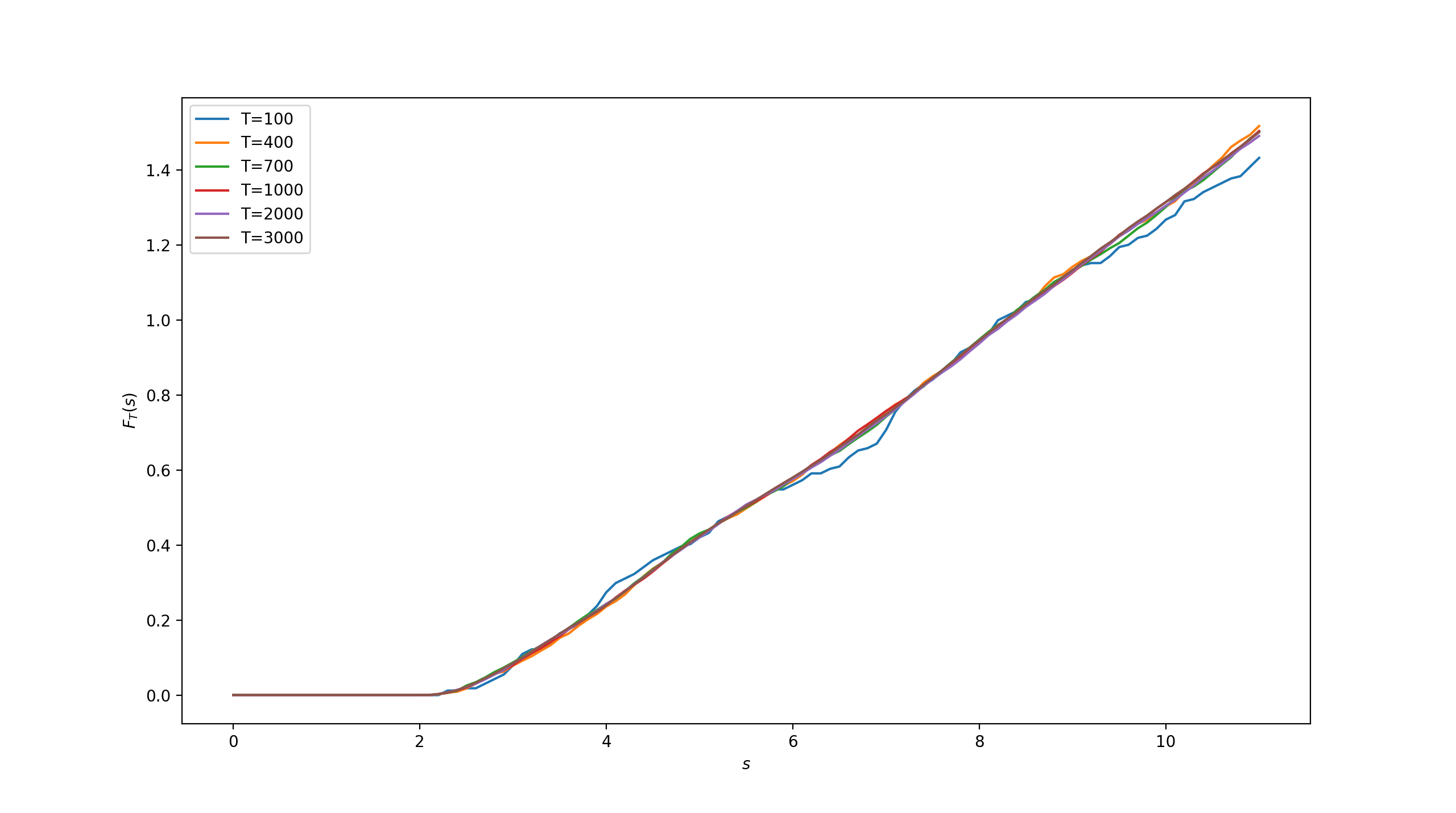}}{\caption{The plot for $F_T$ with various T's}\label{pair}}
             \ffigbox{\includegraphics[scale = 0.27]{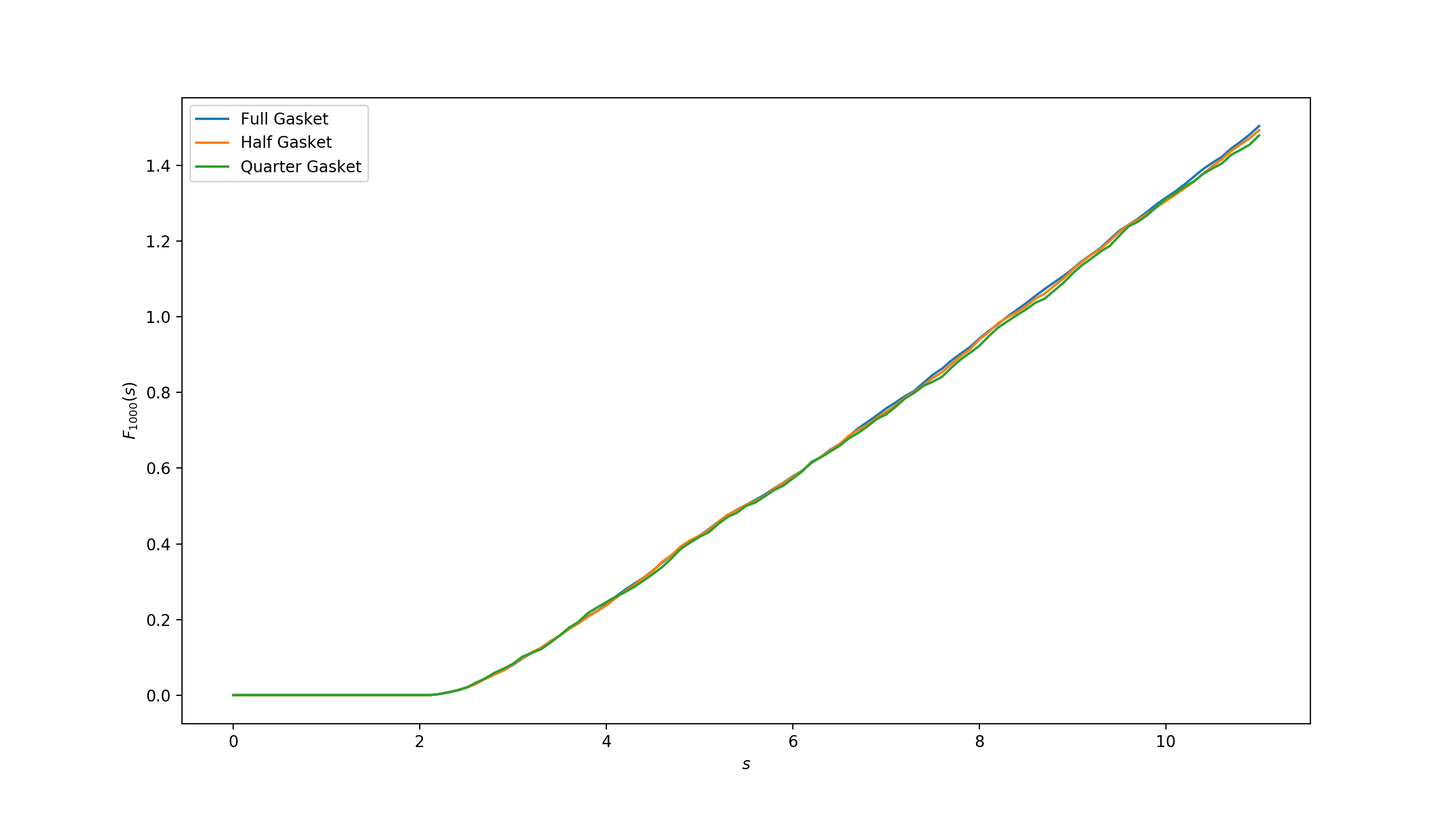}}{\caption{Pair correlation for the whole plane, half plane and the first quadrant}\label{wholehalf}}
           \end{floatrow}
 \end{figure}

Figure \ref{pair} is the empirical plot for the pair correlation function $F_T$ for the gasket $\mathcal P (\frac{1.8}{3}\pi, \frac{3.7}{3}\pi)$, with various $T$ taken. It seems that these curves indeed stay close to each other. \par
We can also study the pair correlation function for centers restricted to some subset $\mathcal{R}$ of $\mathbb C$:
\begin{align}\label{f_t_r}
F_{T,\mathcal R}(s):=\frac{1}{2\#\mathcal C_{T,\mathcal{R}}}\sum_{\substack{p,q\in\mathcal C_{T,\mathcal{R}}\\p\neq q}}\bd{1}\left\{d(p,q)T< s\right\}
\end{align}
where $\mathcal C_{T,\mathcal R}=\mathcal C_T\cap \mathcal R$. \par

By convention if $\mathcal{R}=\mathbb C$, we can also omit $\mathcal R$. Figure \ref{wholehalf} is the plot for $F_{1000,\mathcal R}$ for the gasket $\mathcal P (\frac{1.8}{3}\pi, \frac{3.7}{3}\pi)$, with $\mathcal{R}=\CC$, $\mathcal R=\{z\in\mathbb C\vert \Re z>0\}$ and $\mathcal R=\{z\in\mathbb C\vert \Re z>0,\Im z>0\}$ respectively. The three obtained curves indeed seem close to each other.  \par

\begin{figure}[!h]
           \begin{floatrow}
             \ffigbox{\includegraphics[scale = 0.27]{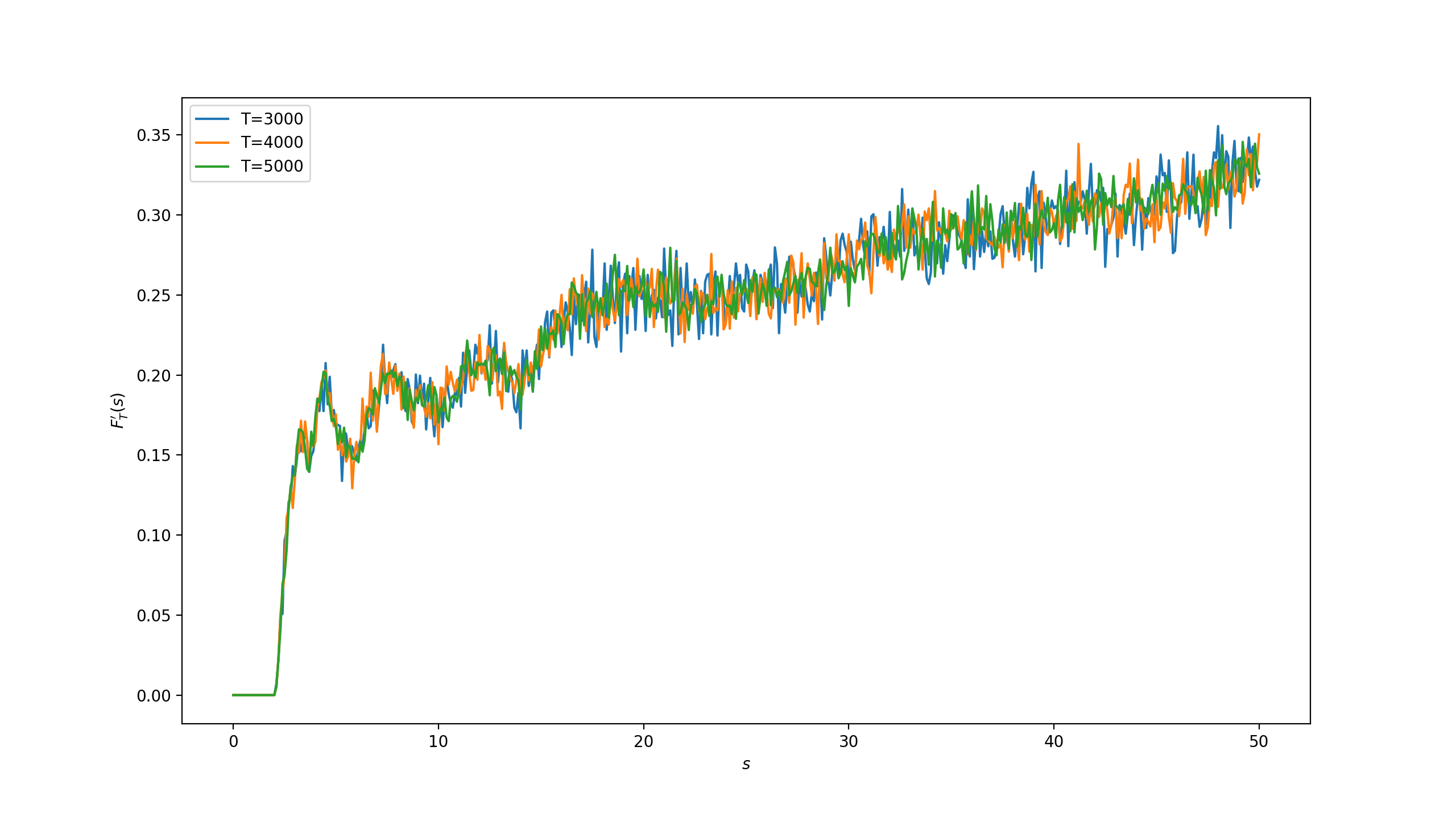}}{\caption{The empirical derivative $F_T'(s)$, with different $T$ taken}\label{deri}}
             \ffigbox{\includegraphics[scale = 0.27]{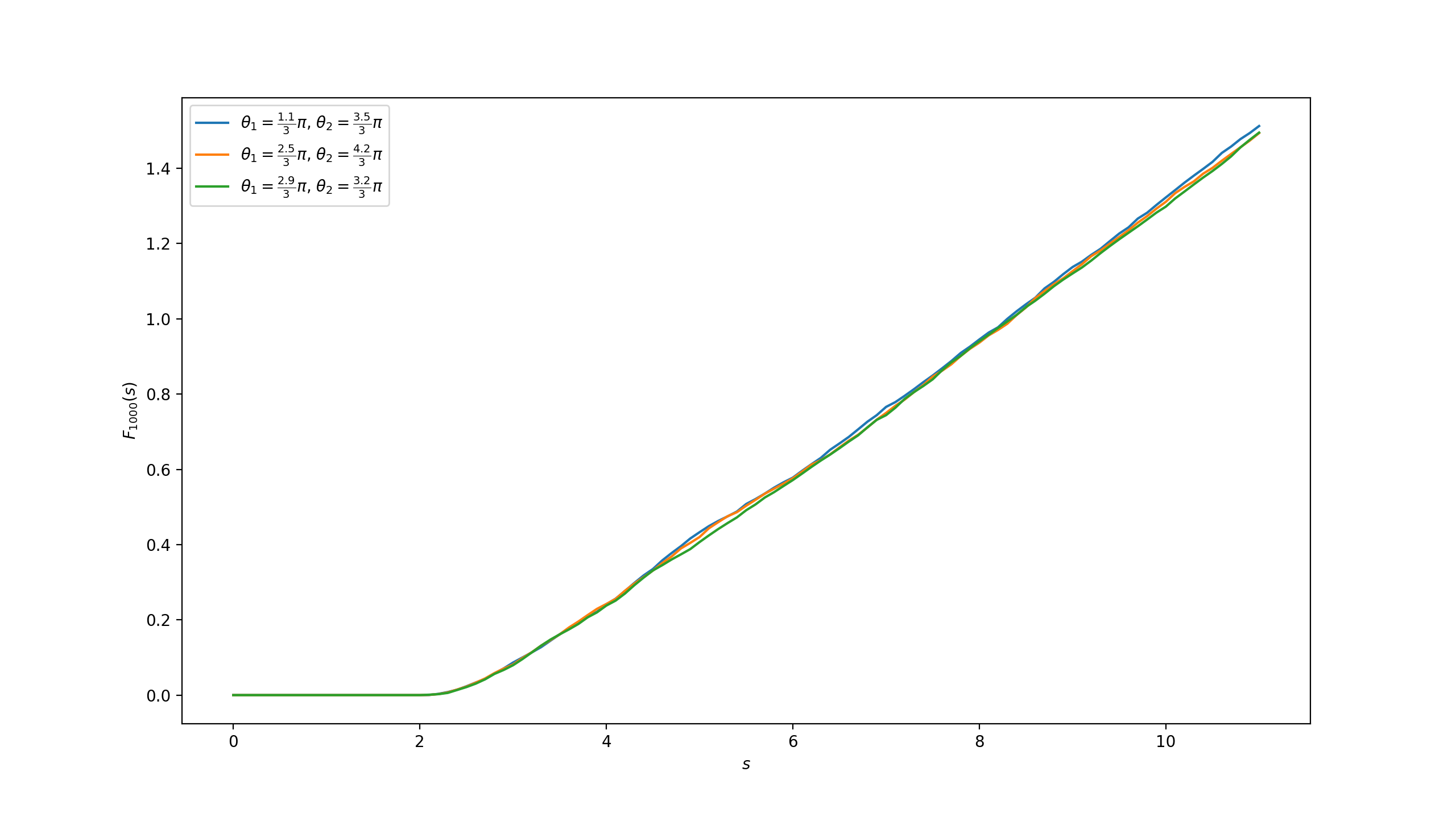}}{\caption{Pair correlation functions for different Apollonian gaskets}\label{figure_1}}
           \end{floatrow}
 \end{figure}

In Figure \ref{deri} we also plot $``{F_T'}(s)"$, the empirical derivative for $F_T(s)$, defined by ${F_T'}(s)=\frac{F_T(s+0.1)-F_T(s)}{0.1}$, for the gasket $\mathcal P (\frac{1.8}{3}\pi, \frac{3.7}{3}\pi)$. Our plot suggests that the derivative of $F_T$ exists and is continuous. The turbulent manner of the plot indicates that a rigorous proof of this claim might be difficult. \par

Figure \ref{figure_1} is the plot for the pair correlation function $F_T$ for three different Apollonian gaskets $\mathcal P(\frac{1.1}{3}\pi, \frac{3.5}{3}\pi)$, $\mathcal P(\frac{2.5}{3}\pi, \frac{3.5}{4.2}\pi)$, $\mathcal P(\frac{2.9}{3}\pi, \frac{3.2}{3}\pi)$. It appears that their limiting pair correlation should be the same.  
 
Based on these findings, we make the following conjecture.
\begin{conj} \label{con1}For any Apollonian gasket $\mathcal{P}$, and any $\mathcal R\subset \mathbb C$ with $\mu(\mathcal R)>0$, there exists a non-negative, monotone, continuously differentiable function $F$ on $[0,\infty)$ which is supported away from 0 such that 
$$\lim_{T\rightarrow\infty}F_{T,\mathcal R}(s)=F(s)$$
for any $s\in[0,\infty)$, where $F_{T,\mathcal{R}}(s)$ is defined in \eqref{f_t_r}. Moreover, the function $F$ is independent of the chosen Apollonian gasket.  
\end{conj}
\subsection{Electrostatic Energy}
The electrostatic energy function $G(T)$ is defined by 
\al{\label{0625}G(T):=\frac{1}{T^{2\delta}}\sum_{\substack{p,q\in\mathcal C_T\\p\neq q}}\frac{1}{d(p,q)}.}
The definition \eqref{0625} agrees with the definition of electrostatic energy for an array of electrons in physics, with an extra normalizing factor $1/T^{2\delta}$.  

\begin{figure}[!h]
\includegraphics[scale=0.28]{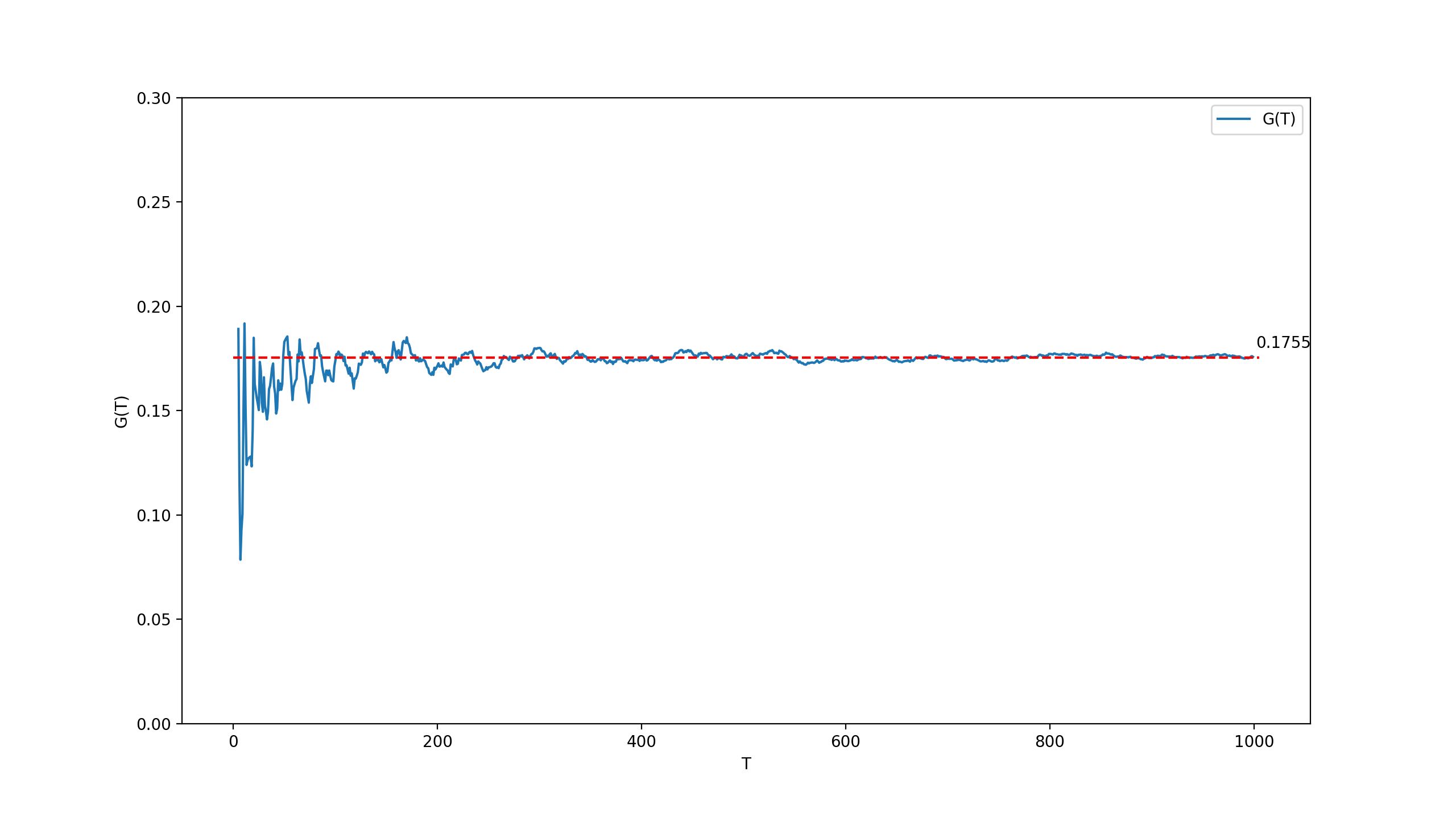}{\caption{The electrostatic energy function $G(T)$ for $\mathcal P(\frac{1.8}{3}\pi, \frac{3.7}{3}\pi)$}\label{elec}}
\end{figure}
Our experiment suggests that $G(T)$ converges to some positive constant when $T$ gets large (see Figure \ref{elec}). We formulate this as a conjecture below. 
\begin{conj}\label{con2}There exists a constant $b>0$, such that 
$$\lim_{T\rightarrow\infty}G(T)=b.$$
\end{conj} 
\subsection{Nearest spacing}

\begin{figure}[!h]
\includegraphics[scale=0.27]{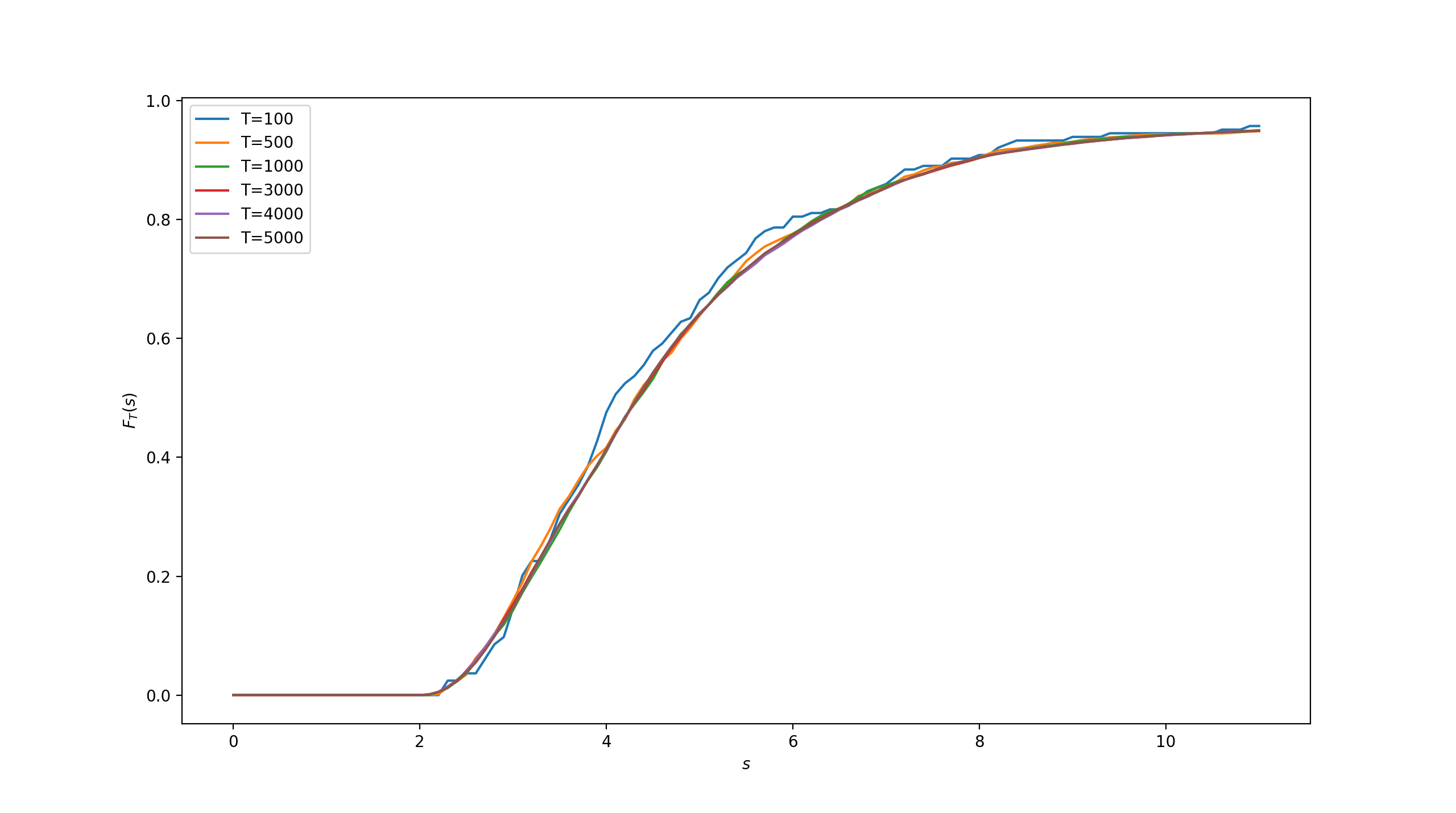}{\caption{The nearest neighbor spacing function $H_T(s)$ for various $T$'s}\label{near}}
\end{figure}

For the set $\mathcal C_T$ and a point $x\in\mathcal C_T$, we let $g_T(x)$ to be the distance between $x$ and a closest point in $\mathcal C_T$ to $x$. The nearest spacing function $H_T(s)$ for the set $\mathcal C_T$ is then defined as 
\begin{align}\label{H_T}
H_T(s):=\frac{1}{\#\mathcal C_T}\sum_{x\in\mathcal C_T}\bd{1}\left\{g_T(x)\cdot T< s\right\}
\end{align}
Again, we need to normalize the distance by multiplying $T$, as we did earlier for the pair correlation function. Figure \ref{near} is the plot of the nearest spacing function for the gasket $\mathcal P(\frac{1.8}{3}\pi, \frac{3.7}{3}\pi)$ for various $T$'s, which also depicts some convergence behavior.   

\begin{conj}\label{con3}There exists a non-negative, monotone, continuous function $H$ on $[0,\infty)$ which is supported away from 0 such that 
$$\lim_{T\rightarrow \infty}H_T(s)=H(s),$$
for any $s\in[0,\infty)$, where $H_T(s)$ is defined as in \eqref{H_T}.
\end{conj}

\section{Conclusion}
Our investigation shows that the spatial statistics of Apollonian gaskets exhibit quite regular behavior, this is probably due to the fact that these gaskets are highly self-symmetric. A possible approach to the proposed conjectures might be via homogeneous dynamics on infinite volume hyperbolic spaces.  \par
\begin{figure}[!h]
           \begin{floatrow}
            
             \ffigbox{\includegraphics[scale = 0.297]{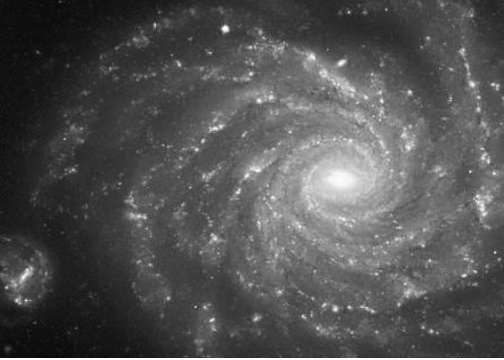}}{\caption{A spiral galaxy}\label{spiral2}}
             \ffigbox{\includegraphics[scale = 0.3]{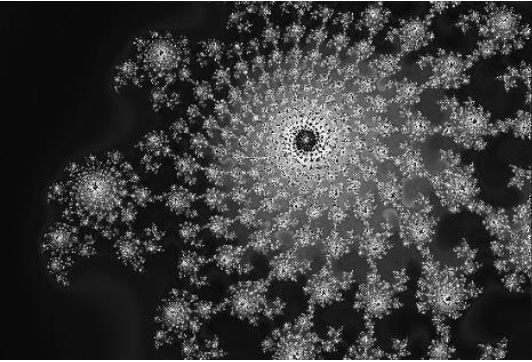}}{\caption{A Mandelbrot set}\label{spiral1}}
           \end{floatrow}
 \end{figure}

  
There are other natural problems on the fine structures of fractal sets. For instance, Figure \ref{spiral2} is the famous Grand Spiral Galaxy (NGC 1232), which can be simulated by a Mandelbrot set constructed from complex dynamics (see Figure \ref{spiral1}). Both pictures are from \cite{Ga13}.\par
 We pose the following question:
 \begin{que}What can one say about the fine structures of stars distribution in a spiral galaxy?
 \end{que}
 \noindent$\bd{Acknowledgement:}$ Our work is the product of an Illinois Geometry Lab (IGL) undergraduate research project during the Spring, 2017 semester, which is supported by the National Science Foundation Grant DMS-1449269. The first three authors are the undergraduate members of our IGL team, the fourth one is a graduate mentor, and the last author is the faculty mentor for the project. We are thankful to Professor Alex Kontorovich for Figures \ref{fig:taxicab} and \ref{apon}, and Professor Gardi for Figures \ref{spiral2} and \ref{spiral1}.

\bibliographystyle{plain}
\bibliography{Apexp}

\end{document}